# New problems on old solitaire boards

George I. Bell and John D. Beasley

```
      - - -
      o o -
      o - o
- - o - o - o - -
- o - o - o - o -
- - o - o - o - -
      o - o
      - o -
      - - -     o - o
                o o o
                o o o
            o o o o o o o
            o o o o o o o
            o o o o o o o
                o o o
                o o o
                o o o       o
                          o o o
                        o o - o o
                      o o o o o o o
                    o + o o o o o o o
                      o o o o o o o
                        o o o o o
                          o o o
                            o
```

# New problems on old solitaire boards

George I. Bell and John D. Beasley

**Abstract**

Some old solitaire boards are brought down from the literature, dusted off, and re-examined, and some remarkable problems are displayed on them.

**Contents**

## 1. Background and history

Solitaire ("Peg Solitaire" in America) is a very familiar one-person board game. The board consists of an array of holes or squares; a man (usually a peg or a marble) is placed in each hole; one man is removed; the rule of play is to jump a man over a neighbouring man, removing the man jumped over; the object is either to reduce to a single man or to leave the men in some specified pattern. But these simple rules yield a game of remarkable subtlety, which has spawned a substantial literature both practical and mathematical. In particular, J. H. Conway published an extensive treatment in 1982 which gave the most important discoveries up to that time (in some of which he had played a major part), and one of us devoted a complete book to the subject a few years later (Berlekamp, Conway and Guy 1982/2004, Beasley 1985/1992). Any statement not more specifically referenced in what follows will be found in one of these works.

The game appears to have originated in France in the late seventeenth century, and to have been the "Rubik's Cube" of the court of Louis XIV (Berey 1697, Trouvain 1698). On the evidence of a passing reference in a letter from Horace Walpole, "Has Miss Harriet found out any more ways at *solitaire*?", it was already established in England in 1746; in 1985, one of us took a very cautious view of this, suspecting that the reference might be to a card game, but David Parlett has written that these fears were groundless: "Patience dates from the late eighteenth century, did not reach England until the nineteenth, and was not called Solitaire when it did" (Walpole 1746, Parlett 1999). Sadly, the widely-quoted legend that the game was invented by a prisoner in the Bastille is almost certainly false. In 1985, the earliest reference to this that could be found was in an English book of 1801, more than a century after the alleged event and in a different country to boot, and nobody has yet brought an earlier reference to our attention.

The game was originally played on the 37-hole board shown in Fig. 1, and a selection of problems on this board was published by Berey (Berey undated). However, it is impossible on this board to solve the puzzle "start by vacating the central hole, play to leave a single man in this hole", and the game is now more usually played on the 33-hole board shown in Fig. 2. But boards of many other shapes and sizes have been tried, and in this paper we shall consider three of these: (a) the 45-hole board of Fig. 3, briefly studied by Johann Christoph Wiegleb in 1779 but little used since; (b) the 39-hole board of Fig. 4, which forms a halfway house between Wiegleb's board and the normal 33-hole board and on which there is a simply stated problem of remarkable difficulty; (c) the 41-hole diamond board of Fig. 5, which received attention in France in the late 19th century (Wiegleb 1779, Lucas 1882/1891). Play from a single initial vacancy to a single survivor is very much harder on these boards than on the 33-hole and 37-hole boards.

Fig. 1: The classical 37-hole board

```
    o o o
  o o o o o
o o o o o o o
o o o o o o o
o o o o o o o
  o o o o o
    o o o
```

Fig. 2: The 33-hole board

```
    o o o
    o o o
o o o o o o o
o o o o o o o
o o o o o o o
    o o o
    o o o
```

Fig. 3: Wiegleb's 45-hole board

```
      o o o
      o o o
      o o o
o o o o o o o o o
o o o o o o o o o
o o o o o o o o o
      o o o
      o o o
      o o o
```

Fig. 4: The 39-hole "semi-Wiegleb" board

```
    o o o
    o o o
    o o o
o o o o o o o
o o o o o o o
o o o o o o o
    o o o
    o o o
    o o o
```

Fig. 5: The 41-hole diamond board

```
        o
      o o o
    o o o o o
  o o o o o o o
o o o o o o o o o
  o o o o o o o
    o o o o o
      o o o
        o
```

Most of the results that follow owe their discovery to computer search, but those that do not involve counting moves have been verified by mathematical analysis. One of us (GIB) programmed the computer to search for solutions of certain kinds, "failure to find" being treated as a provisional but strong indication that no such solution existed; the other (JDB) attacked the unsolved problems analytically, and demonstrated that solutions of the kind required were indeed impossible. A specimen demonstration appears in an appendix. Claims that a certain solution is the shortest possible, or that only a certain number of similarly short solutions exist, are completely computer-dependent and have not been independently verified.

The unsolvability of the problem "start by vacating the central hole, play to leave a single man in this hole" on the 37-hole board is a consequence of a property known as "position class": the various positions possible on a Solitaire board can be divided into 16 different classes, and it is impossible to play from a position in one class to a position in another. The earliest demonstration of this known to us was due to Suremain de Missery and was reported as late as 1842 (Vallot 1842), but the result, once suspected, is not difficult to prove, and we believe that it had in fact been established very early in the history of the game. (It is perhaps proved most simply by marking off the diagonals of the board in threes, mark-mark-clear-mark-mark-clear and so on, and observing that if the total number of pegs in the marked diagonals starts odd it must remain odd throughout the solution and if it starts even it must remain even.) This theory is spelt out in detail in all mathematical treatments of the game, but we shall not go further into it here. All the problems we shall consider will have initial and target positions within the same class, and the reasons for any unsolvability will lie deeper.

**2. Long sweeps on Wiegleb's board 45-hole board**

On the 33-hole board, the longest sweep geometrically possible is a 16-sweep (see for example Fig. 6). However, such a board position cannot be reached from a single vacancy start, or in fact from any starting position with fewer than 16 vacancies (the position itself)! From a single-vacancy start, the best that can be achieved is an 11-sweep, but it is not then possible to play on and reduce to a single man. The longest sweep that can occur as the last move in play from a single-vacancy start to a single-survivor finish has length 9.

Fig. 6: A 16-sweep on the 33-hole board

```
  a b c d e f g
1     o o -
2     o - o
3 - o - o - o -
4 o - o - o - o
5 - o - o - o -
6     o - o
7     - o -
```

Figs. 7a and 7b: The 16-sweeps possible on Wiegleb's board

```
  a b c d e f g h i
1       - - -
2       o o -
3       o - o
4 - - o - o - o - -
5 - o - o - o - o -
6 - - o - o - o - -
7       o - o
8       - o -
9       - - -
```

```
  a b c d e f g h i
1       - - -
2       - o -
3       o - o
4 - - o o o - o - -
5 - o - o - o - o -
6 - - o - o - o - -
7       o - o
8       - o -
9       - - -
```

In these and subsequent figures, "o" denotes a man and "-" an empty hole.

The longest sweep geometrically possible on Wiegleb's board is the same 16-sweep as is available on the 33-hole board. The sweep can begin either from d2 as shown in Fig. 7a or from d4 as in Fig. 7b (or from 10 other locations symmetrically equivalent). Remarkably, each of these 16-loops can be realized as the final move in single-vacancy to single-survivor play.

The easiest way to solve a long-sweep problem is to set up the reverse of the target position, putting a peg where the target position has a hole and vice versa, and to attempt to reduce this position to a single peg. If we succeed, listing the jumps in reverse order gives a solution to the original problem. For example, Fig. 8a shows the reverse of Fig. 7a, and if we set up this position we find (probably after a certain amount of trial and error) that we can play f1-f3, i4-g4-e4, d4-f4-f2, i6-g6-e6-e8, d9-d7, d6-d8, f9-d9-d7, a6-c6, f8-d8-d6-b6, a4-a6-c6-c4, b4-d4, d1-f1-f3-d3-d5-f5-h5, and i5-g5, reducing to a single man at g5. If we now set up a full board, vacate g5, and play i5-g5, f5-h5, d5-f5, and so on, we eventually find ourselves at Fig. 7a, and we can play the spectacular 16-sweep and finish with a single survivor at d2. The reader may care to set up Fig. 8b, which shows the reverse of Fig. 7b, and play to reduce to a single man at g4, though this is appreciably more difficult.

Figs. 8a and 8b: The reverses of Figs. 7a and 7b

```
  a b c d e f g h i       a b c d e f g h i
1       o o o           1       o o o
2       - - o           2       o - o
3       - o -           3       - o -
4 o o - o - o - o o     4 o o - - - o - o o
5 o - o - o - o - o     5 o - o - o - o - o
6 o o - o - o - o o     6 o o - o - o - o o
7       - o -           7       - o -
8       o - o           8       o - o
9       o o o           9       o o o
```

By this means, it can be shown that Fig. 7a can be reached from a full board with just g5 vacated, and Fig. 7b from a full board with just g4 vacated (or d7, by symmetry). Exhaustive search by computer, subsequently confirmed by mathematical analysis, has shown that no other starting positions are possible. The actual solutions discovered by this means are somewhat untidy, since a 6-sweep such as d1-f1-f3-d3-d5-f5-h5 becomes six separate moves f5-h5, d5-f5, d3-d5, f3-d3, f1-f3, and d1-f1 when the jumps are made in reverse order, but once a solution has been found it can easily be rearranged to reduce the number of separate moves. Solutions optimized in this sense appear in Appendix A.

The 16-loop is the longest sweep ending at d2 or d4, but what about the other locations on this board? Although sweeps as long as 13 are geometrically possible, the longest finishing sweep at any other location only has length 9 (to the problem “vacate d9, finish at d3”). It is however possible to have longer *internal* sweeps, and in particular the penultimate move in a solution to the problems “vacate d6 or g6, finish at d3” can be the same 16-sweep from d2 as we have just seen above.

### 3. Optimally short solutions on Wiegleb’s board

Having solved a solitaire problem, the natural next step is to try and minimize the number of separate moves (counting a sequence of jumps by the same man as a single move). On the 33-hole board, a remarkable set of optimal solutions was developed by Ernest Bergholt and Harry O. Davis between 1912 and 1967. One of us demonstrated by computer in 1985 that these solutions were indeed optimal, and this has been confirmed by Jean-Charles Meyrignac. On the 37-hole board, the best pre-computer work was done by Alain Maye, Leonard J. Gordon, and Davis. Computer analysis by Meyrignac subsequently beat four of their solutions by a single move, and demonstrated the remainder to be optimal (Meyrignac 2002).

On Wiegleb’s board, the problem does not appear to have received attention until now. It follows from the “position class” theory that on this board, as on the standard 33-hole board, the starting and finishing holes of a single-vacancy single-survivor problem must be a multiple of three rows and columns apart (for example, if we start by vacating d1, we can hope to finish at d1 itself, a4, d4, g4, or d7, but nowhere else). It follows that there are 36 essentially different single-vacancy single-survivor problems; any other such problem can be transformed into one of these 36 by rotation or reflection.

Wiegleb himself only gave a solution to the problem “vacate d1, finish at d4”, with a note that the inverse problem “vacate d4, finish at d1” was also solvable. In fact each of the 36 problems is solvable except for “vacate e1, play to finish at e1”. One of us established this in 1985 and indicated how the unsolvability of the outstanding case had been proved, but a fully written-out proof would be lengthy and to the best of our knowledge none has been published. Analysis by computer has now confirmed this unsolvability, and has also found the shortest solution to each of the remaining 35 problems. The problem “vacate and finish at e2” requires 23 moves (an interesting parallel with the 33-hole board, where the equivalent problem “vacate and finish at d2” takes at least one move more than any other), and each of the rest can be solved in between 20 and 22 moves. While 16 (nearly half) of the problems can be solved in 20 moves, no 19-move solution has been found.

Table 1 summarizes the results, and specimen solutions are given in Appendix A. Since these results come from lengthy and complex computer runs, they do not constitute a formal proof that the lengths given are the minimum possible. It is possible that some programming bug is present, and as yet the results await independent verification. However, the program has reproduced all the results previously established as shortest possible on the 33-hole board, and we are confident that its results on other boards are similarly correct.

These solutions are very difficult to find by hand, some virtually impossible (though the solutions to “vacate b5 or e5, finish at e2” with their intermediate 12-sweeps were so found by Alain Maye). Peg solitaire problems do not have unique solutions, except in very special cases or on small boards. If two successive moves do not interfere with each other, they can be executed in the opposite order, producing another solution of the same length. Even when one looks at the set of moves in a solution regardless of order, there are usually many different possibilities, all equally valid. However for two of the problems in Table 1 the set of solution moves is in fact unique (to within symmetry, of course), and in three other cases there are only two or three possibilities.

Table 1: Shortest solutions to single-vacancy single-survivor problems on Wiegleb’s board

| Vacate | Finish at | Number of moves | Number of different solutions | Vacate | Finish at | Number of moves | Number of different solutions |
|---|---|---|---|---|---|---|---|
| d1 | d1 | 22 | n/c | e1 | e1 | Unsolvable | |
| d7 | d1 | 20 | 1 | b4, e4 | e1 | 20 | 3 |
| d4 | d1 | 20 | 30 | e7 | e1 | 22 | 39 |
| g4 | d1 | 20 | 38 | | | | |
| a4 | d1 | 21 | 30 | e2 | e2 | 23 | n/c |
| | | | | e8 | e2 | 22 | 1 |
| d2 | d2 | 20 | 775 | b5, e5 | e2 | 22 | 204 |
| d8 | d2 | 20 | 208 | | | | |
| d5 | d2 | 20 | 983 | e3 | e3 | 22 | 213 |
| g5 | d2 | 21 | 683 | e6 | e3 | 21 | 220 |
| a5 | d2 | 22 | 341 | b6 | e3 | 21 | 244 |
| | | | | e9 | e3 | 22 | 2 |
| d3 | d3 | 20 | 364 | | | | |
| d9 | d3 | 20 | 91 | e4, b4 | e4 | 20 | 191 |
| d6 | d3 | 20 | 3798 | e1 | e4 | 22 | 87 |
| g6 | d3 | 20 | 4845 | e7 | e4 | 22 | n/c |
| a6 | d3 | 21 | 2943 | | | | |
| | | | | e5, e2 | e5 | 22 | 3 |
| d4 | d4 | 20 | 40 | | | | |
| g4 | d4 | 20 | 63 | | | | |
| d1 | d4 | 21 | 60 | | | | |

The numbers of solutions are given to within symmetry and ordering of moves; “n/c” indicates “not calculated”.

Specimen solutions are given in Appendix A, and we draw particular attention to the 22-move solution to the central game (vacate and finish at e5). This is the analogue of Bergholt’s 18-move solution on the standard 33-hole board.

**4. A difficult problem on the 39-hole “semi-Wiegleb” board**

The standard 33-hole board and the 45-hole Wiegleb board are special cases of “generalized cross” boards in which a central 3x3 square has a 3xn arm attached to each side. Each arm is two holes deep on the 33-hole board and three holes deep on Wiegleb’s board. The 39-hole board of Fig. 4 offers a halfway house between these boards. A systematic investigation of generalized cross boards by computer has brought to light an unusually difficult single-vacancy single-survivor problem on this board.

Solitaire players have long known that the problem “vacate and finish at d1” (Fig. 9) is the most difficult single-survivor problem on the 33-hole board, and on Wiegleb’s board the equivalent problem “vacate and finish at e1” (Fig. 10) is unsolvable. On the 39-hole board of Fig. 4, the problem “vacate and finish at d1” (Fig. 11) is solvable but only just; apart from the left-right reflection made possible by symmetry, the jumps that we must make are completely determined. As always in a Solitaire solution, we have flexibility in the order in which we make them, but the actual jumps must be the same. The solution, which has been published though not in print (Beasley 2003), appears in Appendix A, and a demonstration of its uniqueness is given in Appendix B. This is the only single-vacancy single-survivor problem we know, on a board of natural shape and reasonable size, where the jumps of a solution are completely determined in this way.

Fig. 9: Difficult

```
  a b c d e f g
1     o - o
2     o o o
3 o o o o o o o
4 o o o o o o o
5 o o o o o o o
6     o o o
7     o o o
```

Fig. 10: Unsolvable

```
  a b c d e f g h i
1       o - o
2       o o o
3       o o o
4 o o o o o o o o o
5 o o o o o o o o o
6 o o o o o o o o o
7       o o o
8       o o o
9       o o o
```

Fig. 11: Uniquely solvable

```
  a b c d e f g
1     o - o
2     o o o
3     o o o
4 o o o o o o o
5 o o o o o o o
6 o o o o o o o
7     o o o
8     o o o
9     o o o
```

In each case, the task is to leave the final survivor in the hole initially empty.

This is quite a different property from the uniqueness or near uniqueness of some of the solutions in the last section, and is much more fundamental. There, we were talking about uniqueness of *moves*, and they were only unique or nearly unique because we were restricting ourselves to solutions of a certain length. Here, we are talking about the constituent *jumps*, and their uniqueness remains whether we combine them into a 21-move solution as is done in Appendix A, or play them all out separately, or do anything in between. However we solve the problem, if we write out the jumps we have made and tick them off one by one against those in Appendix A, we find we have made either exactly the same jumps or symmetrically equivalent ones.

**5. "La corsaire" on the 41-hole diamond board**

One of the solutions given by Berey on the 37-hole board is entitled "Table de la Corsaire". In our notation, he vacates e1, and then plays e3-e1, g3-e3, f5-f3, f2-f4, g5-g3, d3-f3, g3-e3, e4-e2, e1-e3, b3-d3, b5-b3, d5-b5, d7-d5, c1-c3-c5, a3-c3, d3-b3, b2-b4, a4-c4, c5-c3, a5-c5, d5-b5, b6-b4, c7-c5, and f6-d6. This leaves the pattern shown in Fig. 12, after which the man on d1 sweeps off nine other men and e7-c7 finishes the solution. A similar finish is possible to the problem "vacate e7, finish at c7" on the 37-hole board, and also to the problems "vacate c1 or c7, finish at c7" on the 33-hole board.

Fig. 12: "La Corsaire"

```
  a b c d e f g
1     - o -
2   - - o - -
3 - - o - o - -
4 - o - o - o -
5 - - o - o - -
6   - - o - -
7     - - o
```

Fig. 13: "La Corsaire" on the 41-hole diamond board

```
  a b c d e f g h i
1         -
2       - o -
3     - - o - -
4   - - o - o - -
5 - - o - o - o - -
6   - - o - o - -
7     - - o - -
8       - - o
9         -
```

Fig. 14: A deceptive setting of "La Corsaire"

```
  a b c d e f g h i
1         o
2       o o o
3     o o - o o
4   o o o o o o o
5 o + o o o o o o o
6   o o o o o o o
7     o o o o o
8       o o o
9         o
```

Play on the 41-hole diamond board is vastly more difficult than on the 37-hole board (the only solvable single-vacancy single-survivor problems are "vacate f8 or c5, play to finish at d8 or g5" and problems equivalent to these), and we were therefore surprised to discover that a "corsaire" finish is possible on this board as well. The natural equivalent of Fig. 12 on the 41-hole diamond board is shown in Fig. 13, and if we start by vacating f8 or c5 we can indeed play to this position.

There is more. If we look at the possible moves of the man at e2 earlier in the play, we find that it can start at c4, g4, c6, g6, or e8, but not at e2 itself. The same is of course true of other orientations of the problem; for example, if we vacate e3 (it's always nice to have the initial vacancy on the vertical axis of symmetry) and play for a corsaire finish h5-...-b5 and b6-b4, we find that the corsaire man can start at f3, f7, d3, d7, or b5, but not at h5 itself. We can therefore pose the problem in the deceptive form shown in Fig. 14: "Vacate e3, mark the man at b5 in some way, and play to reduce to a single survivor, this marked man making a 9-sweep at the penultimate move".

A moderately advanced player, who has read or worked out the "position class" theory and knows that an initial vacancy at e3 means that any survivor must finish at b4 or h4 (or e7, which is not relevant here), will try to set up a finale b5-...-h5 and h6-h4, and he will not succeed; no such finale is possible. Instead, the marked man must migrate to h5 earlier in the play, and then perform a corsaire sweep back to b5. A solution with this property is given in Appendix A.

## 6. Summary

The game of Peg Solitaire may have a history of over three hundred years, but it is very far from exhausted. This paper has presented some of the more interesting problems that have recently come to light. We trust it has given pleasure, and we hope it may prompt others to search for the further delights that are surely still awaiting discovery.

## Appendix A: Solutions

Although most of the solutions which follow were originally generated by computer, the moves of some have been reordered to give a more natural progression round the board. A much more extensive set of Solitaire solutions, on these and on other boards, can be found on the web site <http://www.geocities.com/gibell.geo/pegsolitaire/>.

### Long-sweep solutions on Wiegleb's board

Vacate g5, finish at d2 with a 16-sweep, solution derived from Fig. 8a rearranged to minimize the number of separate moves: i5-g5, f5-h5, d5-f5, d3-d5, f3-d3, f1-f3, f4-f2, d1-f1-f3 (8), b4-d4-f4, c6-c4, a6-c6, d6-b6, a4-a6-c6 (13), d8-d6, f8-d8, e6-e8, d9-d7, d6-d8, f9-d9-d7 (19), g6-e6, i6-g6, g4-e4, i4-g4 (23), d2-f2-f4-h4-h6-f6-f8-d8-d6-d4-f4-f6-d6-b6-b4-d4-d2.

Vacate g4, finish at d4 with a 16-sweep. Fig. 8b can be reduced to a single man at g4 by playing a6-c6, d6-b6, a4-a6-c6-c4, b4-d4, f9-f7, f6-f8, d9-f9-f7-d7, d1-d3-d5, d8-d6-d4, i4-g4-e4-e2, f2-d2, f1-d1-d3-d5-f5-h5, h6-h4, i6-i4-g4. The key to this solution lies in the parallel moves g4-e4-e2 and d3-d5-f5-h5, and any solution must contain either these moves or equivalent ones across the SE corner (i4-g4-e4-e6-e8, a6-c6, d6-b6, a4-a6-c6-c4, d1-d3, b4-d4-d2, f1-d1-d3-f3, f9-f7-f5, f2-f4-f6, d8-f8, d9-f9-f7-f5-h5 etc). The jumps can be rearranged to give the following 22-move solution to the original problem: i4-g4, h6-h4, f5-h5, d5-f5, d3-d5, d6-d4, d1-d3-d5, f2-d2, f1-d1-d3 (9), b4-d4, c6-c4, a6-c6, d8-d6-b6, a4-a6-c6 (14), f7-d7, f9-f7, f6-f8, d9-f9-f7 (18), e4-e2, g4-e4, i6-i4-g4 (21), d4-d2-f2-f4-h4-h6-f6-f8-d8-d6-b6-b4-d4-f4-f6-d6-d4.

These solutions, in 24 and 22 moves respectively, are the shortest possible.

The problem "Vacate d9, finish at d3 with a 9-sweep" can be solved in a simple systematic way by playing f9-d9, e7-e9, e5-e7, g5-e5, f7-f5, e5-g5, h6-f6, h4-h6, i6-g6, f6-h6, f4-h4, i4-i6-g6-g4, f2-f4, d2-f2, f1-f3, f4-f2, d4-d2, d1-d3, h4-f4-d4-d2, b4-d4, c6-c4, d4-b4, a4-c4, d6-d4-b4, a6-a4-c4, b6-b4-d4, d8-d6, d9-f9-f7-d7-d5-d3-d1-f1-f3-d3. A rearrangement saving a move is possible (instead of h4-f4-d4-d2, b4-d4, c6-c4, d4-b4, play b4-d4-d2, c6-c4, h4-f4-d4-b4).

The problem "Vacate d6 or g6, finish at d3 with an internal 16-sweep" can be solved most simply by starting from Fig. 8a and playing a6-c6-e6-e4-e2, f2-d2, d1-d3-d5-b5 (again these parallel moves), b4-b6, f9-f7, f6-f8, d9-f9-f7-d7, d8-d6, i4-g4, f4-h4, i6-i4-g4-g6, h6-f6, and either a4-a6-c6-e6 and f6-d6 or a4-a6-c6-e6-g6. This leaves two pegs untouched at e1 and f1 but otherwise reduces the board to a single peg at d6 or g6. If we now refill the board, vacate g6 or d6, and play these jumps in reverse order, we come down to Fig. 7a with two extra pegs at e1 and f1, and we can complete the solution by playing the 16-sweep and the move f1-d1-d3. A rearrangement in 23 moves is possible.

### Optimally short solutions on Wiegleb's board

We arrange the solutions in order of finishing hole, as in Table 1.

*Finish at d1*

Vacate d1: d3-d1, f2-d2, f4-f2, f1-f3 (4), d4-f4, e6-e4, c6-e6, a6-c6 (8), d8-d6-b6, f6-d6, c4-c6-e6, f7-d7, f9-f7, h6-f6-f8, d9-f9-f7 (15), g4-g6, i4-g4, f4-h4, i6-i4-g4 (19), a4-a6-c6, b4-b6-d6-d8-f8-f6-h6-h4-f4-d4-d6-f6-f4-f2, d1-f1-f3-d3-d1 (22).

Vacate d4 or g4: f4-d4 or e4-g4, then c4-e4, a4-c4 (3), f6-f4-d4-b4, d2-d4, f2-d2, h4-f4-f2, f1-f3 (8), g6-g4, i6-g6, f8-f6-h6, i4-i6-g6 (12), d8-f8, f9-f7, d6-f6-f8, d9-f9-f7 (16), a6-a4-c4-e4-e6-e8, c6-c4, b6-b4-d4-d6-d8-f8-f6-h6-h4-f4-f2, d1-f1-f3-d3-d1 (20).

Vacate d7: d9-d7, d6-d8, d4-d6, f4-d4, c4-e4, f2-f4-d4, a4-c4-e4 (7), b6-b4, a6-a4-c4, f9-d9-d7-d5-b5 (10), d2-f2, f1-f3 (12), h4-f4-f2, g6-g4, i6-g6, f6-d6-b6-b4-d4-f4-f6-h6, i4-i6-g6 (17), e8-e6-e4, f8-f6-h6-h4-f4-d4-d2, d1-f1-f3-d3-d1 (20). As a 20-move solution, this is unique to within move ordering.

Vacate a4: c4-a4, b6-b4, a4-c4 (3), d5-b5, d7-d5, d9-d7 (6), d4-d6-d8, d2-d4-b4, a6-a4-c4 (9), f4-d4-b4-b6-d6, f2-d2,

h4-f4-f2, f1-f3 (13), h6-h4, i4-g4, f6-f4-h4, i6-i4-g4 (17), f9-d9-d7-d5-f5-h5, f7-d7, f8-d8-d6-f6-h6-h4-f4-f2, d1-f1-f3-d3-d1 (21).

*Finish at d2*

Vacate d2 or d5: d4-d2 or d3-d5, then d1-d3 (2), f4-d4-d2, e6-e4, c6-e6, a6-c6, d8-d6-b6, f6-d6, c4-c6-e6, a4-a6-c6 (10), f7-d7, f9-f7, h6-f6-f8, d9-f9-f7 (14), g4-g6, i4-g4, f2-f4-h4, i6-i4-g4 (18), f1-d1-d3-f3, b4-b6-d6-d4-f4-f6-d6-d8-f8-f6-h6-h4-f4-f2-d2 (20).

Vacate d8: d6-d8, d9-d7 (2), d4-d6-d8, b4-d4, c6-c4, a6-c6, f6-d6-b6, a4-a6-c6 (8), d3-d5, d1-d3 (10), e8-e6, h6-f6-d6-d4-d2, g4-g6, i4-g4, f4-h4, i6-i4-g4 (16), f2-f4-h4-h6-f6-f4-d4-b4-b6-d6, f9-d9-d7-d5-f5, f1-d1-d3-f3, f8-f6-f4-f2-d2 (20).

Vacate g5: e5-g5, e3-e5, c4-e4, a4-c4 (4), f3-f5, f1-f3, h4-f4-f2, d1-f1-f3 (8), g6-g4, i6-g6, f6-h6, i4-i6-g6 (12), f8-f6-f4-d4-b4, d7-f7, d9-d7, d6-d8, f9-d9-d7 (17), b6-d6-d4, e6-e4-c4-c6, a6-a4-c4 (20), and finish with the 12-loop d2-f2-f4-h4-h6-f6-f8-d8-d6-b6-b4-d4-d2 (21).

Vacate a5: c5-a5, e5-c5, d3-d5, d1-d3, d6-d4-d2, f1-d1-d3 (6), f6-d6, c6-e6, a6-c6 (9), d8-d6-b6, f7-d7, f9-f7, h6-f6-f8, d9-f9-f7 (14), g4-g6, i4-g4, f4-h4, i6-i4-g4 (18), f2-d2-d4-f4-f6-d6-d8-f8-f6-h6-h4-f4-f2, b4-d4, a4-a6-c6-c4-e4-e2, f2-d2 (22).

*Finish at d3*

Vacate d3: d5-d3, d7-d5, d9-d7 (3), b4-d4-d6-d8, c6-c4, a6-c6, f6-d6-b6, a4-a6-c6 (8), f4-d4, e2-e4-e6, h4-f4-f6-d6, g6-g4, i6-g6, f8-f6-h6, i4-i6-g6 (15), d3-d5-d7-f7, d1-d3, f9-d9-d7, f2-f4-h4-h6-f6-f8-d8-d6-b6-b4-d4-d2, f1-d1-d3 (20).

Vacate d9: d7-d9, d5-d7, f8-d8-d6, f6-f8, f9-f7, f4-f6-f8, d9-f9-f7-d7-d5-f5 (7), d3-d5, b4-d4, c6-c4, a6-c6, h6-f6-d6-b6, a4-a6-c6 (13), g4-g6, i4-g4, f2-f4-h4, i6-i4-g4 (17), d1-d3-f3, d4-b4-b6-d6-d4-f4-f6-h6-h4-f4-f2-d2, f1-d1-d3 (20).

Vacate d6 or g6: f6-d6 or e6-g6, then c6-e6, h6-f6-d6, c4-c6-e6, a6-c6, d8-d6-b6, a4-a6-c6 (7), f7-d7, f9-f7, f4-f6-f8, d9-f9-f7 (11), g4-g6, i4-g4, f2-f4-h4, i6-i4-g4 (15), d5-f5, d3-d5, d1-d3-f3 (18), b4-b6-d6-d8-f8-f6-h6-h4-f4-f6-d6-d4-f4-f2-d2, f1-d1-d3 (20). These solutions include an internal 14-sweep.

Vacate a6: c6-a6, e6-c6, d8-d6-b6, d4-d6, f4-d4, a6-c6-e6-e4, b4-b6, a4-a6-c6 (8), f7-d7, f9-f7, f6-f8, d9-f9-f7 (12), h6-f6, g4-g6, i4-g4, f2-f4-h4, i6-i4-g4 (17), d3-d5-b5, d1-d3-f3, f6-h6-h4-f4-f6-f8-d8-d6-b6-b4-d4-f4-f2-d2, f1-d1-d3 (21).

*Finish at d4*

Vacate d4 or g4: f4-d4 or e4-g4, then c4-e4, a4-c4 (3), h4-f4-d4-b4, d2-d4, f3-d3, f1-f3, f6-f4-f2, d1-f1-f3 (9), g6-g4, i6-g6, f8-f6-h6, i4-i6-g6 (13), d8-f8, f9-f7, d6-f6-f8, d9-f9-f7 (13), a6-a4-c4-e4-e6-e8, c6-c4, b6-b4-d4-d2-f2-f4-h4-h6-f6-f8-d8-d6-d4 (20).

Vacate d1: d3-d1, f2-d2, d1-d3, d4-d2, f1-d1-d3 (5), e4-e2, g4-e4, i4-g4 (8), b4-d4-d2-f2-f4-h4, d6-d4-f4, c6-c4, a6-c6, d8-d6-b6, a4-a6-c6 (14), f8-d8, d9-d7, f6-d6-d8, f9-d9-d7 (18), i6-i4-g4-e4-e6-e8, g6-g4, h6-h4-f4-f6-f8-d8-d6-b6-b4-d4 (21).

*Finish at e1*

Vacate b4 or e4: d4-b4 or c4-e4, then d6-d4, d3-d5, d1-d3 (4), b6-d6-d4-d2, a4-c4-c6, d8-d6-b6, f6-d6, a6-c6-e6 (9), f7-d7, f9-f7, h6-f6-f8, d9-f9-f7 (13), g4-g6, i4-g4, f4-h4, i6-i4-g4 (17), f2-f4-f6-f8-d8-d6-f6-h6-h4-f4, f1-d1-d3-f3-f5-d5, a5-c5-e5-e3-e1 (20). The only alternatives still allowing a 20-move ordering are to play f6-d6-b6 and d8-d6 at moves 7 and 8, or f7-d7, f9-f7, f6-f8, h6-f6-d6, a6-c6-e6 at moves 8-12.

Vacate e7: e5-e7, c6-e6, a6-c6, e3-e5, c4-e4, a4-c4 (6), d8-d6-d4, d3-d5, d1-d3 (9), f6-d6-d4-d2, f7-d7, f9-f7, h6-f6-f8, d9-f9-f7 (14), g4-g6, i4-g4, f4-h4, i6-i4-g4 (18), f2-f4-h4-h6-f6-f8-d8-d6-b6-b4-d4, f5-d5-b5, f1-d1-d3-d5, a5-c5-e5-e3-e1 (22).

*Finish at e2*

Vacate e2: e4-e2, c4-e4, a4-c4 (3), d2-d4, d5-d3, f2-d2-d4-b4, f4-f2, f1-f3, h4-f4-f2, f6-f4, d1-f1-f3-f5 (11), g6-g4, i6-g6, f8-f6-h6, i4-i6-g6 (15), d7-f7, d9-d7, d6-d8, f9-d9-d7 (19), e5-e3, b5-d5, b6-d6-d8-f8-f6-h6-h4-f4-f6-d6-d4, a6-a4-c4-e4-e2 (23).

Vacate e8: e6-e8, c6-e6, d4-d6, d7-d5, d9-d7, b4-d4-d6-d8, f9-d9-d7 (7), a6-c6-c4 (8), d2-d4-b4, f3-d3, f1-f3, f4-f2, d1-f1-f3 (13), h4-f4-f2-d2-d4, g6-g4, i6-g6, f6-h6, i4-i6-g6 (18), f8-d8-d6-f6-f4-h4-h6-f6, f7-f5-d5, a5-c5-e5-e3, a4-c4-e4-e2 (22). As a 22-move solution, this is unique to within symmetry and move ordering.

Vacate b5 or e5: d5-b5 or c5-e5, then d3-d5, f4-d4, f6-f4, f3-f5, f1-f3, h4-f4-f2, d1-f1-f3 (8), g6-g4, i6-g6, f8-f6-h6, i4-i6-g6 (12), d7-f7, d9-d7, d6-d8, f9-d9-d7 (16), c4-e4, a4-c4, d2-f2-f4-f6-f8-d8-d6-f6-h6-h4-f4-d4-b4, a6-a4-c4 (20), b6-b4-d4-d6, c6-e6-e4-e2 (22). This solution, with its internal 12-sweep, was found by Alain Maye.

*Finish at e3*

Vacate e3: e5-e3, c4-e4, d6-d4, f6-d6, c6-e6, a6-c6 (6), d3-d5, d8-d6-d4, f7-d7, f9-f7, h6-f6-f8, d9-f9-f7 (12), g4-g6, i4-g4, f4-h4, i6-i4-g4 (16), f2-f4-f6-d6-b6, a4-a6-c6 (18), d1-d3-f3, b4-b6-d6-d8-f8-f6-h6-h4-f4-f2-d2, f1-d1-d3-d5, c5-e5-e3 (22).

Vacate b6 or e6: d6-b6 or c6-e6, then f6-d6, a6-c6-e6, c4-c6, e5-c5 (5), then join the previous solution after move 6 (21 moves in all).

Vacate e9: e7-e9, e5-e7, c6-e6, a6-c6 (4), d4-d6, d7-d5, d9-d7, d2-d4-d6-d8, b4-b6-d6, f9-d9-d7-d5-b5 (10), f3-d3, f1-f3, f4-f2, d1-f1-f3 (14), h4-f4-f2-d2-d4-b4-b6, g6-g4, i6-g6, f6-h6, i4-i6-g6 (19), f8-f6-f4-h4-h6-f6-d6, a4-a6-c6-e6, e7-e5-e3 (22). The only alternative still allowing a 22-move ordering is to play b4-b6-d6-d8 and d2-d4-d6 at moves 8 and 9.

*Finish at e4*

Vacate e4 or b4: c4-e4 or d4-b4, then a4-c4 (2), f4-d4-b4, f6-f4, f3-f5, f1-f3 (6), d6-f6-f4-f2, h4-f4, g6-g4, i6-g6, f8-f6-h6, i4-i6-g6 (12), d7-f7, d9-d7, d2-d4-d6-d8, f9-d9-d7 (16), b5-d5-f5-f3-d3, d1-f1-f3, b6-d6-d8-f8-f6-h6-h4-f4-f2-d2-d4, a6-a4-c4-e4 (20).

Vacate e1: e3-e1, e5-e3, c4-e4, d6-d4, d3-d5, d1-d3 (6), f4-d4-d2, f6-d6-d4, h6-f6-f4, h4-h6, f3-f5-h5, f1-d1-d3-f3 (12), i6-g6, f8-f6-h6, i4-i6-g6 (15), d7-f7, d9-d7, b6-d6-d8, f9-d9-d7 (19), a4-c4-c6, f2-f4-h4-h6-f6-f8-d8-d6-b6-b4, a6-a4-c4-e4 (22).

Vacate e7: Follow moves 1-20 of "Vacate g5, finish at d2" above, rotated 90 degrees clockwise, and finish with h4-h6-f6-f8-d8-d6-b6-b4-d4-d2-f2-f4, g4-e4 (22).

*Finish at e5*

Vacate e5 or e2: e3-e5 or e4-e2, then c4-e4, d6-d4, d3-d5, d1-d3, b6-d6-d4-d2, f1-d1-d3 (7), a4-c4-c6, d8-d6-b6, f6-d6, a6-c6-e6 (11), f7-d7, f9-f7, h6-f6-f8, d9-f9-f7 (15), g4-g6, i4-g4, f4-h4, i6-i4-g4 (19), f2-d2-d4-f4-h4-h6-f6-f8-d8-d6-f6-f4, f3-f5-d5, a5-c5-e5 (22). The only alternatives still allowing a 22-move ordering are to play f6-d6-b6 and d8-d6 at moves 9 and 10, or f7-d7, f9-f7, f6-f8, h6-f6-d6, a6-c6-e6 at moves 10-14. Bergholt's 18-move solution to the equivalent problem "vacate d4 or d1, finish at d4" on the 33-hole board is d2-d4 or d3-d1, f3-d3, e1-e3, e4-e2, e6-e4 (5), g5-e5, d5-f5, g3-g5-e5, c3-e3, a3-c3, b5-d5-f5-f3-d3-b3 (11), c1-e1-e3-e5, c7-c5, c4-c6 (14), e7-c7-c5, a5-a3-c3, c2-c4-c6-e6-e4-c4, b4-d4.

**39-hole "semi-Wiegleb" board**

Vacate and finish at d1. It is shown in Appendix B that we need the following jumps: d3-d1 twice, d5-d3 twice, d8-d6, d4-d2, a4-a6, a6-c6 twice, d6-b6, b4-d4, g4-g6, g6-e6 twice, d6-f6, f4-d4, c5-e5 (or e5-c5, everything that follows being reflected left to right), b5-d5, f5-d5, c1-e1, e1-e3 twice, e4-e2, e6-e4, c2-c4, c7-c5, e9-c9, c9-c7 twice, c6-c8, c4-c6, c7-c5 again, e3-e5 twice, e6-e4 twice more, e8-e6. An optimal ordering is given by d3-d1, d5-d3, f4-d4-d2, e6-e4, e3-e5, e1-e3, g6-e6 (7), b5-d5, c7-c5, c9-c7, a6-c6 (11), b4-d4, d6-b6, c2-c4-c6-c8, d8-d6-f6, e8-e6-e4-e2 (16), a4-a6-c6, e9-c9-c7-c5-e5, g4-g6-e6-e4, c1-e1-e3-e5, f5-d5-d3-d1 (21).

**41-hole diamond board**

The problem "Vacate f8 or c5, play to the position of Fig. 13" can be solved most simply by playing d8-f8, e6-e8, c7-e7 or c7-c5, e6-c6, d8-d6, then g6-e6, g4-g6, e4-g4, c4-e4, c6-c4, e6-c6 (it can be shown that any solution must contain this cycle of six jumps, either this way round or in the reverse direction, and it is simplest to play them straight away), then e8-e6, g7-e7, e6-e8, e9-e7, h6-f6, i5-g5, f2-f4, d3-f3, e1-e3, f3-d3, f5-f3, h4-f4, g3-e3, e4-e2, d2-d4, d5-d3, c3-e3, b4-d4, a5-c5, b6-d6. Rearrangements in 24 moves are possible, giving 26 moves to reduce to a single survivor. One of them, rotated through 90 degrees, appears in the next paragraph.

The problem shown in Fig. 14 can be solved by playing c3-e3, d5-d3, d2-d4, b5-d5-d3 (4), b4-d4, c7-c5, d7-d5-b5, a5-c5 (8), f7-d7, d8-d6, e9-e7 (11), f5-f7, h6-f6, g4-g6, i5-g5 (15), e4-g4, f2-f4, d3-f3-f5-h5 (18), e1-e3, h4-f4, g7-g5, e6-g6-g4, g3-g5, f8-f6 (24). We now have Fig. 13 rotated through 90 degrees and the man initially at b5 is now at h5, and we can finish by playing h5-f5-d5-d7-f7-f5-f3-d3-d5-b5 and b6-b4 (26). This solution is also an optimally shortest solution to the problem "vacate e3, finish at b4"; indeed, each of the single-vacancy single-survivor problems solvable on this board has a solution in 26 moves but none shorter.

## Appendix B: A specimen non-computer analysis

In principle, any unsolvable Solitaire problem can be so proved by trying every possibility in turn and verifying that none works, and this is just what a computer search does. But this is rarely practicable by hand, and non-computer analyses normally use techniques developed in the 1960s by J. M. Boardman, J. H. Conway, and R. L. Hutchings (Berlekamp, Conway and Guy 1982/2004, Beasley 1985/1992). As a specimen, we show the uniqueness (to within symmetry and order of jumps) of the solution to "vacate and finish at d1" on the 39-hole "semi-Wiegleb" board.

Our basic approach will be to write down the numbers in Fig. 15 (–1 at d1, +1 everywhere else), and to try to apply successive adjustments "–1, –1, +1" to the numbers in adjacent holes so that they are all eventually reduced to zero. In effect, we shall play Solitaire with numbers in a table rather than with pegs on a board. The advantage of this approach is that once we see that a particular jump will be necessary, for example because it is the only way to bring a negative number up to zero, we can apply it to the table straight away; we do not need to wait until a suitable configuration of pegs arises on the board.

Fig. 15: The task in numerical form

| | a | b | c | d | e | f | g |
|---|---|---|---|---|---|---|---|
| 1 | | | 1 | -1 | 1 | | |
| 2 | | | 1 | 1 | 1 | | |
| 3 | | | 1 | 1 | 1 | | |
| 4 | 1 | 1 | 1 | 1 | 1 | 1 | 1 |
| 5 | 1 | 1 | 1 | 1 | 1 | 1 | 1 |
| 6 | 1 | 1 | 1 | 1 | 1 | 1 | 1 |
| 7 | | | 1 | 1 | 1 | | |
| 8 | | | 1 | 1 | 1 | | |
| 9 | | | 1 | 1 | 1 | | |

We now note that the first jump of a solution must be d3-d1, and the next jump d5-d3; also that the last jump will be d3-d1 again, and the penultimate jump d5-d3. If we apply these jumps to the table, we get Fig. 16:

Fig. 16: After considering d3-d1 twice and d5-d3 twice

| | a | b | c | d | e | f | g |
|---|---|---|---|---|---|---|---|
| 1 | | | 1 | 1 | 1 | | |
| 2 | | | 1 | -1 | 1 | | |
| 3 | | | 1 | 1 | 1 | | |
| 4 | 1 | 1 | 1 | -1 | 1 | 1 | 1 |
| 5 | 1 | 1 | 1 | -1 | 1 | 1 | 1 |
| 6 | 1 | 1 | 1 | 1 | 1 | 1 | 1 |
| 7 | | | 1 | 1 | 1 | | |
| 8 | | | 1 | 1 | 1 | | |
| 9 | | | 1 | 1 | 1 | | |

Next, we assign a value to each hole on the board as shown in Fig. 17. We shall be adding up the values of the holes which are occupied, and these values have the property that if A, B, C are any three adjacent holes in line and f(A) etc are their values then f(A)+f(B) is at least as great as f(C). Since the effect of a jump from A over B into C is to replace the contribution f(A)+f(B) by a contribution f(C), *the sum of the values of the holes occupied can never increase.*

Fig. 17: An assignment of values to holes

| | a | b | c | d | e | f | g |
|---|---|---|---|---|---|---|---|
| 1 | | | 0 | 0 | 0 | | |
| 2 | | | 0 | 1 | 0 | | |
| 3 | | | 0 | 0 | 0 | | |
| 4 | -1 | 1 | 0 | 1 | 0 | 1 | -1 |
| 5 | 0 | 0 | 0 | 0 | 0 | 0 | 0 |
| 6 | -1 | 1 | 0 | 1 | 0 | 1 | -1 |
| 7 | | | 0 | 0 | 0 | | |
| 8 | | | 0 | 1 | 0 | | |
| 9 | | | 0 | 0 | 0 | | |

If we evaluate the task shown by Fig. 16 according to the values in Fig. 17, we find we have contributions –1 from d2/a4/d4/g4/a6/d6 and +1 from b4/g4/b6/d6/f6/d8, total zero. But the value of our target position (all zeros) is also zero, and it follows that we can never make a jump which reduces our evaluation according to Fig. 17; once this evaluation has become negative, we can never get it back up again. In particular, a jump over d8 will have just this effect (it will replace a contribution 0+1 by 0), so our solution cannot contain such a jump. But we must clear d8 somehow, and the only remaining candidate is d8-d6. Apply this to Fig. 16, and also the jump d4-d2 which is needed to bring d2 up to zero, and we have Fig. 18.

Fig. 18: The task after considering d8-d6 and d4-d2

| | a | b | c | d | e | f | g |
|---|---|---|---|---|---|---|---|
| 1 | | | 1 | 1 | 1 | | |
| 2 | | | 1 | 0 | 1 | | |
| 3 | | | 1 | 0 | 1 | | |
| 4 | 1 | 1 | 1 | -2 | 1 | 1 | 1 |
| 5 | 1 | 1 | 1 | -1 | 1 | 1 | 1 |
| 6 | 1 | 1 | 1 | 2 | 1 | 1 | 1 |
| 7 | | | 1 | 0 | 1 | | |
| 8 | | | 1 | 0 | 1 | | |
| 9 | | | 1 | 1 | 1 | | |

Fig. 19 shows a second assignment of values to holes. Again, f(A)+f(B) is at least as great as f(C) whenever A, B, C are any three adjacent holes in line, so once more the sum of the holes occupied can never increase. But this time, if we evaluate the current task (shown by Fig. 18) according to the values assigned, we find we have a net total of +1 (there are contributions –4 from d4, –2 from a4/g4, –1 from d5/a6/g6, +2 from b4/f4/d6, and +1 from a5/b5/f5/g5/b6/f6). In other words, we can afford to lose 1 (and indeed we must lose it at some stage, since we need eventually to reduce everything to zero), but we cannot afford to lose 2.

Fig. 19: A second assignment of values to holes

| | a | b | c | d | e | f | g |
|---|---|---|---|---|---|---|---|
| 1 | | | 0 | 0 | 0 | | |
| 2 | | | 0 | 2 | 0 | | |
| 3 | | | 0 | 0 | 0 | | |
| 4 | -2 | 2 | 0 | 2 | 0 | 2 | -2 |
| 5 | 1 | 1 | 0 | 1 | 0 | 1 | 1 |
| 6 | -1 | 1 | 0 | 1 | 0 | 1 | -1 |
| 7 | | | 0 | 0 | 0 | | |
| 8 | | | 0 | 1 | 0 | | |
| 9 | | | 0 | 0 | 0 | | |

Now, reverting to Fig. 18, how are we going to clear a5? The only candidate jumps are a4-a6, a6-a4, and a5-c5, but each of the two latter loses 2 when we evaluate according to Fig. 19 and we have just seen that we cannot afford this. So we must play a4-a6, and by similar arguments we must play a6-c6 twice, d6-b6 (the two jumps a6-c6 have left us with –1 at b6, and to play b4-b6 would lose 2 according to Fig. 19), and b4-d4. We shall need equivalent jumps g4-g6, g6-e6 twice, d6-f6, and f4-d4 on the right-hand side as well, and if we apply all these jumps to Fig. 18 we get Fig. 20:

Fig. 20: The task after considering a4-a6 etc and g4-g6 etc

| | a | b | c | d | e | f | g |
|---|---|---|---|---|---|---|---|
| 1 | | | 1 | 1 | 1 | | |
| 2 | | | 1 | 0 | 1 | | |
| 3 | | | 1 | 0 | 1 | | |
| 4 | 0 | 0 | 0 | 0 | 0 | 0 | 0 |
| 5 | 0 | 1 | 1 | -1 | 1 | 1 | 0 |
| 6 | 0 | 0 | 2 | 0 | 2 | 0 | 0 |
| 7 | | | 1 | 0 | 1 | | |
| 8 | | | 1 | 0 | 1 | | |
| 9 | | | 1 | 1 | 1 | | |

We still have value 1 according to Fig. 19 (+1 at b5/f5, –1 at d5). We need to reduce this to zero, and the only way of doing this without incurring an intolerable loss according to Fig. 17 is to play c5-e5 or the symmetrically equivalent e5-c5. For present purposes, let us suppose c5-e5. We shall now need b5-d5 and f5-d5 to bring the number in d5 up to zero, and we have Fig. 21:

Fig. 21: After considering c5-e5, b5-d5, and f5-d5

| | a | b | c | d | e | f | g |
|---|---|---|---|---|---|---|---|
| 1 | | | 1 | 1 | 1 | | |
| 2 | | | 1 | 0 | 1 | | |
| 3 | | | 1 | 0 | 1 | | |
| 4 | 0 | 0 | 0 | 0 | 0 | 0 | 0 |
| 5 | 0 | 0 | -1 | 0 | 1 | 0 | 0 |
| 6 | 0 | 0 | 2 | 0 | 2 | 0 | 0 |
| 7 | | | 1 | 0 | 1 | | |
| 8 | | | 1 | 0 | 1 | | |
| 9 | | | 1 | 1 | 1 | | |

We now have non-zero numbers only at d1/d9 and on the c and e files, and we can proceed a little faster. We need to clear d1, and the only candidates are e1-c1 and c1-e1. Try e1-c1. There follows c1-c3 twice, c4-c2 to liquidate the resulting deficiency at c2, c6-c4 to liquidate the deficiency at c4, c7-c5 twice to liquidate the double deficiency at c5... No, it's not possible.

All right, try c1-e1. We must now play e1-e3 twice, e4-e2, e6-e4, c2-c4, and c7-c5, and we have reduced to Fig. 22:

Fig. 22: After considering c1-e1 etc

| | a | b | c | d | e | f | g |
|---|---|---|---|---|---|---|---|
| 1 | | | 0 | 0 | 0 | | |
| 2 | | | 0 | 0 | 0 | | |
| 3 | | | 0 | 0 | 2 | | |
| 4 | 0 | 0 | 1 | 0 | 0 | 0 | 0 |
| 5 | 0 | 0 | 0 | 0 | 0 | 0 | 0 |
| 6 | 0 | 0 | 1 | 0 | 1 | 0 | 0 |
| 7 | | | 0 | 0 | 1 | | |
| 8 | | | 1 | 0 | 1 | | |
| 9 | | | 1 | 1 | 1 | | |

To clear d9, try c9-e9; no, the surpluses on c4/c6/c8 prove intractable. Try e9-c9; c9-c7 twice, c6-c8, c4-c6, c7-c5, and we have Fig. 23:

Fig. 23: After considering e9-c9 etc

| | a | b | c | d | e | f | g |
|---|---|---|---|---|---|---|---|
| 1 | | | 0 | 0 | 0 | | |
| 2 | | | 0 | 0 | 0 | | |
| 3 | | | 0 | 0 | 2 | | |
| 4 | 0 | 0 | 0 | 0 | 0 | 0 | 0 |
| 5 | 0 | 0 | 0 | 0 | 0 | 0 | 0 |
| 6 | 0 | 0 | 0 | 0 | 1 | 0 | 0 |
| 7 | | | 0 | 0 | 1 | | |
| 8 | | | 0 | 0 | 1 | | |
| 9 | | | 0 | 0 | 0 | | |

Now, at last, we have a choice:

(a) a pair of jumps e2-e4/e4-e2 across e3, e6-e8, a pair e9-e7/e7-e9 across e8;

(b) a pair e2-e4/e4-e2 across e3, e8-e6, a pair e7-e5/e5-e7 across e6;

(c) a pair e2-e4/e4-e2 across e3, e8-e6, a pair f6-d6/d6-f6 across e6;

(d) e3-e5 twice, e6-e4 twice, e8-e6.

However, this has been playing Solitaire with numbers. If we revert to the real board, we find that the region e1/e2/e3 starts full and finishes empty, so we need two jumps outwards across its boundary (one to broach it initially, one to remove the last man from it), and only option (d) provides them.

**Acknowledgments**

Our grateful thanks to Jean-Charles Meyrignac for many computational tips and discussions, and to Alain Maye, who is able to solve many of these problems without the aid of a computer.

George I. Bell
5040 Ingersoll Pl.
Boulder
CO 80303
USA

gibell@comcast.net

John D. Beasley
7 St James Road
Harpenden
Herts AL5 4NX
UK

johnbeasley@mail.com